\documentstyle[12pt]{article}
\newcommand{\const}{\mathop{\rm const}\limits}

\newcommand{\Var}{\mathop{\rm Var}\limits}

\newcommand{\cov}{\mathop{\rm cov}\limits}

\newcommand{\vraisup}{\mathop{\rm vraisup}\limits}

\newcommand{\mod}{\mathop{\rm mod}\limits}

\newcommand{\Ent}{\mathop{\rm Ent}\limits}

\begin{document}

\begin{center}

{\bf  CLT for continuous random processes \\
under approximations terms. } \par

\vspace{4mm}

 $ {\bf E.Ostrovsky^a, \ \ L.Sirota^b } $ \\

\vspace{4mm}

$ ^a $ Corresponding Author. Department of Mathematics and computer science, Bar-Ilan University, 84105, Ramat Gan, Israel.\\
\vspace{3mm}
E - mail: \ galo@list.ru \  eugostrovsky@list.ru\\

\vspace{4mm}

$ ^b $  Department of Mathematics and computer science. Bar-Ilan University,
84105, Ramat Gan, Israel.\\

\vspace{3mm}

E - mail: \ sirota3@bezeqint.net\\

\vspace{5mm}

 {\sc Abstract.} \\

\end{center}

\vspace{4mm}

  We formulate and prove a new sufficient conditions for Central Limit Theorem  (CLT) in the space
of continuous functions in the terms typical for the approximation theory.\par
 We prove that the {\it sufficient} conditions for continuous CLT obtained by N.C.Jain and M.B.Marcus are under
some natural additional conditions {\it necessary}. \par
 We provide also some examples in order to show the exactness of obtained results and illustrate briefly the applications
 in the Monte-Carlo method. \par

  \vspace{3mm}


\vspace{4mm}

 \par

\vspace{4mm}

\section{Notations. Statement of problem.}

\vspace{3mm}

  Let $ \eta(t), \ t \in T = [0, 2 \pi]   $ be separable centered (zero mean) random process (r.p.) with finite covariation
 function  $ R(t,s) = \cov(\eta(t), \eta(s)) = {\bf E} \eta(t) \eta(s),  \  \{\eta_i(t) \}, \ i = 1,2,\ldots $  be a
sequence of independent copies (identical distributed) of $ \eta(t). $ We will denote as usually the probabilistic notions
by $ {\bf P}, \ {\bf E}, \ \cov, \ \Var $ etc. \par
 Let as introduce the following sequence of a random processes:

 $$
 \zeta_n(t) = n^{-1/2} \sum_{j=1}^n \eta_j(t). \eqno(1.1)
 $$
 Let also $ C(T) $ be a Banach space of all periodical: $ f(0) = f( 2 \pi)  $ continuous  functions equipped with
usually uniform norm   $  ||f|| = \max_{t \in T} |f(t)|, \  \zeta_{\infty}(t) = \zeta(t)  $ be a centered separable Gaussian
process with at the same covariation  function $  R(t,s). $ It is clear that the finite-dimensional distributions of r.p.
$ \zeta_n(t) $ converges as $ n \to \infty $ to the ones for $ \zeta_{\infty}(t). $\par

\vspace{3mm}

{\bf Definition 1.1.} {\it We will say as ordinary that the r.p. } $ \eta(t) $ {\it satisfies the  Central Limit Theorem in the
space } $ C(T), $  {\it notation:}   $  \eta(\cdot) \in CLT(C(T)),  $
{\it if for any continuous bounded functional } $ F: C(T) \to R $

$$
\lim_{n \to \infty} {\bf E} F(\zeta_n(\cdot)) =  {\bf E} F(\zeta_{\infty}(\cdot)).
$$

\vspace{3mm}

 {\bf  Our purpose  in this  article is to formulate and prove  sufficient conditions for CLT in the space $  C(T) $ under
the terms  which are habitual in the theory of approximation. }\par

\vspace{3mm}

 This problem has a long history; see e.g. \cite{Dudley1}, \cite{Kozachenko1}, \cite{Marcus2},
\cite{Ledoux1}, \cite{Ostrovsky1}, \cite{Ostrovsky106}. About applications CLT in
the space $  C(T)  $ in the Monte-Carlo method and in statistics see, e.g. \cite{Frolov1},  \cite{Grigorjeva1}, \cite{Ostrovsky106}.\par
 Indeed, if  $  \eta(\cdot) \in CLT(C(T)),  $ then for all positive values $  u $

 $$
\lim_{n \to \infty} {\bf P} \left( \max_{t \in T} |\zeta_n(t)| > u \right)  =
{\bf P} \left( \max_{t \in T} |\zeta_{\infty}(t)| > u \right) .
 $$

 This problem contains as a subproblem the finding of a sufficient condition for continuity with probability one of random processes,
Gaussian or not, see \cite{Bednorz1} - \cite{Bednorz4}, \cite{Fernique1} - \cite{Fernique3},   \cite{Garsia1},  \cite{Heinkel1},
\cite{Ibragimov1}, \cite{Kwapien1}, \cite{Marcus2}- \cite{Marcus1}, \cite{Ostrovsky100} - \cite{Ostrovsky104},
\cite{Ostrovsky207} - \cite{Ostrovsky210}, \cite{Ostrovsky301}, \cite{Ral'chenko1}, \cite{Talagrand1} - \cite{Talagrand5},
\cite{Watanabe1}.\par
  Note   that in the articles \cite{Dudley1}, \cite{Kozachenko1}, \cite{Marcus2}, \cite{Ostrovsky1} the CLT is formulated in the
so-called {\it metric entropy } terms,  in  \cite{Dudley1}, \cite{Heinkel1}, \cite{Ostrovsky106}, \cite{Ostrovsky301} by means of the notions
{\it  majorizing measure}.\par

\vspace{3mm}
 The paper is organized as follows. In the next section
 we formulate and prove a new sufficient conditions for Central Limit Theorem  (CLT) in the space
of continuous functions. In the third section we prove that the {\it sufficient} conditions for continuous CLT
obtained by N.C.Jain and M.B.Marcus are still under some natural additional conditions {\it necessary}. \par
 We provide also some examples in the $  4^{th} $ section in order to show the exactness of obtained results and to illustrate
briefly the applications in the Monte-Carlo method.  The last section contains some concluding remarks.\par

\vspace{4mm}
 We need to introduce some useful notations.   We will use the well-known  Vallee-Poussin sums, which play a very important role in
the approximation theory, see e.g. \cite{Poussin1}, \cite{Timan1}, chapter 5. \par

 Recall that the Vallee-Poussin kernel $ K_{n,p}(t) $ is defined as follows:
$$
K_{n,p}(t) = \frac{\sin((2n+1-p)t/2) \cdot \sin((p+1)t/2)}{2(p+1) \sin^2t/2}, \ p \in [1,2,\ldots,n].
$$
 It is known that $ K_{n,p}(t)  $ is trigonometrical polynomial of degree  $ n: \ K_{n,p}(t) \in A(n), $
where $ A(n) $  denotes the set of all $ 2 \pi $ periodical trigonometrical polynomials of degree $ \le n. $    \par

 The Vallee-Poussin  approximation (sum) $ V_{n,p}[f](t) $ for a periodical integrable function $ f $
may be defined as follows:

$$
V_{n,p}[f](t) := [f * K_{n,p}](t)
$$
(periodical convolution). Recall that

$$
V_{n,p}[f](t) = \frac{1}{p+1} \sum_{k=n-p}^n S_k[f](t),
$$
where $ S_k[f](t)  $ is the $ k^{th} $ partial Fourier sum for the function $ f. $ \par
 We pick hereafter for definiteness for the values $ n \ge 4  \ p = p(n)
\stackrel{def}{=} \Ent(n/2) := [n/2] ,  $ (integer part), so that  $ V_{n,p}[f](t)  \in A(n) $ and

$$
|| f(\cdot) - V_{n,p(n)}[f](\cdot) || \le C \cdot E([n/2],f),
$$
where $ E(m,f) $ denotes the error in the uniform norm of the best trigonometrical approximation
of the function $  f  $ by means of trigonometrical polynomials with degree $ \le m: $

$$
E(m,f)  = \inf_{g \in A(m)} \max_{t \in T}|f(t) - g(t)| =  \inf_{g \in A(m)} ||f(t) - g(t)||,
$$
see \cite{Natanson1}, chapter 6.  \par
  We define for arbitrary r.p. $ \xi(t) $ the so-called {\it generating functional}
 $ \Phi(\xi; \psi)$ as follows:

 $$
\Phi_{\xi} = \Phi(\xi; \psi) = {\bf E} e^{\int_T \xi(s) \ d \psi(s)},  \eqno(1.2)
 $$
if there exists. Here $ \psi(s), \ s \in [0,2 \pi] $  is  any {\it deterministic} function of bounded variation. \par
 Further, we denote by $ \cal{N} $ a set of all strictly increasing  sequences of natural numbers
 $ \{ n(k) \}, \ k=1,2, \ldots, \ n(1) = 1; $ and define  for each such a sequence \\ $ \vec{n} = \{ n(k)\} \in \cal{N} $

 $$
 W^{(k)}(t) = W_{\vec{n}}(n(k+1), n(k))(t) := V_{n(k+1), p(n(k+1))}(t) -  V_{n(k), p(n(k))}(t), \eqno(1.3)
 $$
 and we define for arbitrary periodical random process $ \xi =\xi(t), \ t \in T $

 $$
 Z_{k} [\xi](t)= Z_{ \vec{n},k} [\xi](t)= [W^{(k)}*\xi](t) =
 $$
 $$
 [V_{n(k+1), p(n(k+1))}* \xi](t) -  [V_{n(k), p(n(k))}* \xi](t),\eqno(1.4)
 $$

 $$
\Psi_{\vec{n}}(\xi;\lambda, n(k), n(k+1))= \Psi(\xi;\lambda, n(k), n(k+1)) =
 $$

 $$
  (2 \pi)^{-1} \int_T \  {\bf E} e^{ \lambda Z_k[\xi](t) } \ dt, \ \lambda = \const > 0;
 $$
 evidently, $ \Psi_{ \vec{n}}(\cdot;\cdot, \cdot, \cdot) $  may be easily expressed through the  generating functional $  \Phi(\cdot). $
  Namely, let

 $$
 \psi_t^{(k)}(s) = \int_0^s W^{(k)}(t-x) \ dx,
 $$
 then

 $$
 \Psi(\xi; \lambda, n(k), n(k+1)) = (2 \pi)^{-1} \int_T {\bf E} e^{ \lambda \int_T \xi(s) W^{(k)}(t-s) ds } \ dt =
 $$

 $$
 (2 \pi)^{-1} \int_T dt \ {\bf E} \ e^{ \lambda  \ \int_T \xi(s) \ d \psi_t^{(k)}(s)} =
 (2 \pi)^{-1} \int_T  \ \Phi(\xi;\lambda \psi_t^{(k)}(\cdot) ) \ dt.\eqno(1.5)
 $$

 \vspace{3mm}

  Define also at last
  $$
  U(\Phi_{\xi}; n(k), n(k+1)) = U_{\vec{n}}(\xi; n(k), n(k+1)) =
  $$

  $$
  \inf_{\lambda > 0} \left[  \frac{\log n(k+1)  +  \log \Psi_{\vec{n}}(\xi; \lambda, n(k), n(k+1))}{\lambda} \right].\eqno(1.6)
  $$

\vspace{4mm}

\section{ Main result: sufficient conditions for CLT for continuous processes.}

\vspace{4mm}

{\bf A. Weak compactness of a family of random processes.} \par
\vspace{3mm}

  Let $ \xi_{\alpha}(t), \ \alpha \in \cal{A}  $ be a {\it family} of separable stochastically continuous
periodical processes, $  \cal{A}  $ be arbitrary set. Assume that for some non-random point $ t_0 \in T $ the family
of one-dimensional r.v. $ \xi_{\alpha}(t_0)  $ is tight.  \par

\vspace{4mm}

{\bf Theorem 2.1.}\\
 {\it 1. \ Let} $ \alpha \in \cal{A} $   {\it be a given. If for some sequence }  $  \{n(k) \} \in \cal{N} $

$$
\sum_{k=1}^{\infty} U(\Phi_{\xi_{\alpha}}; n(k), n(k+1) ) < \infty, \eqno(2.1)
$$
{\it then  almost all trajectories of r.p.}  $ \xi_{\alpha}(t) $ {\it are continuous.} \par

\vspace{3mm}

{\it 2. \ Suppose that there exists a single sequence } $ \{ n(k) \}  $
{\it such that for all the set } $ \cal{A} $ {\it the series (2.1) are uniform equiconvergent:  }

$$
\lim_{m \to \infty} \sup_{\alpha \in \cal{A}}  \sum_{k=m}^{\infty} U(\Phi_{\xi_{\alpha}}; n(k), n(k+1) ) = 0. \eqno(2.2)
$$
{\it Then the family of distributions $ \mu_{\alpha}(\cdot) $ in the space $ C(T) $ generated by} $ \xi_{\alpha}(\cdot):  $

$$
 \mu_{\alpha}(A) = {\bf P} (\xi_{\alpha}(\cdot) \in A )
$$
{\it is weakly  compact in this space.} \par

\vspace{3mm}

{\bf Proof.} \\
 Let  $ \alpha \in \cal{A} $ be a fix.
We will use one of the main results of the article \cite{Ostrovsky301}:

$$
{\bf E} || Z_k[\xi_{\alpha}] ||  \le  C \ U(\Phi_{\xi_{\alpha}}; n(k), n(k+1)  ), \eqno(2.3)
$$
where $  C  $ is an absolute constant. \par

 We conclude by virtue of condition (2.1) that  the following series  converges:

$$
\sum_{k=1}^{\infty} {\bf E}  || Z_k[\xi_{\alpha}] ||  < \infty,
$$
with him

$$
\sum_{k=1}^{\infty}  || Z_k [\xi_{\alpha}]  ||  < \infty \ (\mod {\bf P}).
$$
 Therefore, the partial sums  of $ Z_k [\xi_{\alpha}](t), $ namely, the  sequence of the r.p.

$$
 \sum_{m=1}^k  Z_m [\xi_{\alpha}](t)  =  [V_{n(k+1), p(n(k+1))}* \xi](t) \eqno(2.4)
$$
converges uniformly on $  t, \ t \in T $ also with probability one. The limiting as $ k \to \infty $
process coincides with $  \xi_{\alpha}(t) $  since it is continuous in probability; thus, it is
sample part continuous. We proved the  first proposition of theorem 1.1. \par
 Further, let the condition (2.2) be satisfied. We denote

$$
\epsilon(m) = \sup_{\alpha \in \cal{A}}  \sum_{k=m}^{\infty} U(\Phi_{\xi_{\alpha}}; n(k), n(k+1) ), \eqno(2.5)
$$
then  $  \lim_{m \to \infty} \epsilon(m) = 0. $ \par
 As long as

$$
 \sum_{m=1}^k  Z_m [\xi_{\alpha}](t)  =  [V_{n(k+1), p(n(k+1))}* \xi](t)
$$
is  as the function on the variable $ t $   the trigonometrical polynomial  of degree less than
$ n(k+1),  $ we conclude

$$
\sup_{\alpha \in \cal{A}} {\bf E} E(n(k+1),\xi_{\alpha}) \le C_1 \ \epsilon(n(k)) \to 0, \ k \to \infty. \eqno(2.6)
$$

 We can use the {\it inverse}  theorems of approximation theory; see for example,
Stechkin's  estimate (\cite{Timan1}, chapter 6, section 6.1:)

$$
\omega(f, 1/n) \le \frac{C_2}{n} \cdot  \sum_{m=0}^n E_m(f).
$$
 We ensue for the module of continuity $ \omega(\xi_{\alpha}(\cdot), 1/n)$
   for  periodical r.p. $ \xi_{\alpha}(t): $ \par

$$
\sup_{\alpha \in \cal{A}}{\bf E} \omega(\xi_{\alpha}(\cdot), 1/n) \le \frac{C_3}{n} \cdot \sum_{m=0}^n {\bf E}  \ E_m(\xi).
$$
 Note as a  consequence taking into account the monotonicity of the function $ m \to E(m,f) : $

$$
 \lim_{n \to \infty} \sup_{\alpha \in \cal{A}} {\bf E} \ \omega(\xi_{\alpha}(\cdot), 1/n) = 0 \eqno(2.7)
$$
 and therefore

 $$
 \forall \epsilon > 0 \ \Rightarrow  \lim_{\delta \to 0+} \sup_{\alpha \in \cal{A}}{\bf P} ( \omega(\xi_{\alpha}(\cdot), \delta) > \epsilon ) = 0.
 $$
  Our proposition follows now from theorem 1 in the book  of I.I.Gikhman and A.V.Skorohod
 \cite{Gikhman1},  chapter 9, section 2; see also \cite{Prokhorov1},  after applying the Tchebychev's inequality.\par

{\bf Remark 2.1.} The conditions of theorem 2.1 are essentially non-improvable still for the Gaussian processes, see
\cite{Fernique1},  \cite{Marcus1}, \cite{Ostrovsky209}, \cite{Ostrovsky210}, \cite{Ostrovsky211}, \cite{Ostrovsky301}. \par

\vspace{4mm}

{\bf B. CLT for  random processes.} \par
\vspace{3mm}

{\bf Theorem 2.2.} \\

{\it \ Suppose that there exists a single strictly increasing sequence } $ \{ s(k) \}  $
{\it of natural numbers such that the following series (2.8) are uniform equiconvergent:  }

$$
\lim_{m \to \infty} \sup_{n}  \sum_{k=m}^{\infty} U \left(\Phi^n_{\eta, \psi(\cdot)/\sqrt{n}}; s(k), s(k+1) \right) = 0. \eqno(2.8)
$$
 {\it Then  the r.p. } $  \eta(t), \ t \in T $   {\it satisfies CLT in the space C(T)  }.\par
\vspace{3mm}

{\bf Proof.} \ The generating functional  for the r.p. $ \zeta_n(\cdot), $ i.e. $ \Phi_{\zeta_n}(\psi) $
may be calculated and uniformly estimated as follows:

$$
  \Phi_{\zeta_n}(\psi) =  \Phi^n_{\eta}(\psi/\sqrt{n}), \eqno(2.9)
$$

$$
  \Phi_{\zeta_n}(\psi) \le \sup_n \left[  \Phi^n_{\eta}(\psi/\sqrt{n}) \right] < \infty.
$$

  It remains to apply the second assertion of theorem 2.1, choosing $  \{\alpha\} = \{ 1,2,\ldots.  \}. $\par

\vspace{4mm}

\section{Necessary conditions for CLT for continuous processes.}

\vspace{4mm}
 In this section  the set $  T  $ be instead $ [0, 2 \pi]  $ the arbitrary compact metric space equipped with the distance
function $ d = d(s,t). $  The r.p. $ \eta(t) $ remains to be separable and centered.\par

  N.C.Jain and M.B.Marcus in  \cite{Marcus2} are formulated and proved in particular the following
famous result (we retell its in our notations). \par
{\bf Theorem of Jain and Marcus.} Assumptions: there exists a random variable $  M  $ such that

 $$
| \eta(t) - \eta(s)| \le M \cdot \rho(t,s), \eqno(3.1)
 $$
{\it condition of factorization}; where $ \rho(t,s) $ is  continuous
deterministic distance, more exactly, semi-distance, possible different on the source distance $  d $  on the set $ T; $ \\

$$
{\bf E} M^2 < \infty, \eqno(3.2)
$$
{\it moment condition}; \par

$$
\int_0^1 H^{1/2}( T,\rho, z) \ dz < \infty, \eqno(3.3)
$$
{\it entropy condition. } Here $ H(\rho, T,\epsilon)  $ denotes a metric entropy function of the set $ T $ relative the distance
$ \rho(\cdot, \cdot) $  at the point $ \epsilon, \ 0 < \epsilon < 1. $\par

 Then the random process (field) $ \eta(t) $ satisfies the CLT in the space $ C(T,d). $\par

\vspace{3mm}

{\it We discus in this section the necessity of all conditions 3.1; 3.2; 3.3 for the CLT in the space $ C(T,d). $  \par
 So, it will be presumed sometimes further that the r.p. $ \eta(t) $ satisfies this CLT.}\par
\vspace{3mm}
{\bf Proposition 3.1.} Assume that the r.p. $ \eta(t) $ is sample part continuous,   still without  CLT. Then the condition
of factorization (3.1) is satisfied. \par
 This assertion follows immediately from the main result of the articles \cite{Buldygin1},  \cite{Ostrovsky302}; see also
 \cite{Buldygin2}. \par
\vspace{3mm}
{\bf Proposition 3.2.} Suppose in addition $ \eta(t) $  satisfies the CLT in $ C(T,d); $  then in the factorization (3.1)
the metric $ \rho(\cdot,\cdot)  $ may be selected such that

$$
\forall  p \in (0,2) \ \Rightarrow  \ {\bf E} |M|^p < \infty. \eqno(3.4)
$$
 {\bf Proof.} It is proved in the famous book of M.Ledoux and M.Talagrand \cite{Ledoux1}, chapter  10, section 1, p. 274-277
that if the r.p. $ \eta(\cdot) $  satisfies the CLT in the separable Banach space $ X $ with the norm $ ||\cdot||X, $ then
$ {\bf E}||\eta||^p < \infty. $ Therefore

$$
{\bf E} \sup_{t \in T} |\eta(t)|^p < \infty, \ p \in (0,2).
$$
Note that the function $ u \to |u|^p $ is Young - Orlicz function satisfying the well-known  $ \Delta_2 $ condition.\par
 It remains to use on of the main result of aforementioned article \cite{Ostrovsky302}. \par
  Analogously may be proved the following assertion. \\

\vspace{3mm}
{\bf Proposition 3.3.} Assume that the r.p. $ \eta(t) $ is sample part continuous,   still without  CLT. Moreover,
let

$$
{\bf E} \sup_{t \in T} |\eta(t)|^2 < \infty.
$$
 Then the distance $ \rho $  in the factorization representation (3.1) may be selected such that  the r.v. $ M $
has finite second moment:  $ {\bf E} M^2 < \infty.  $\par

\vspace{3mm}

 Before formulating the next result, we introduce together with R.M.Dudley \cite{Dudley1}, N.C.Jain and M.B.Marcus  \cite{Marcus2}
 the third distance $  \tau = \tau(t,s)  $ on the set $  R  $, also natural, as follows:

$$
\tau(t,s) = \sqrt{ \Var(  \eta(t) - \eta(s))} = \left[ {\bf E} (\eta(t) - \eta(s))^2 \right]^{1/2}.\eqno(3.5)
$$
 It follows from the condition (3.1) that $  \tau(t,s) \le C \cdot \rho(t,s). $\par

\vspace{3mm}

{\bf Proposition 3.4.} Suppose the r.p. $ \eta(t), \ t \in [0, 2 \pi] $  satisfies the CLT in $ C(T,d), $
is stationary in wide sense.  Assume also that the metric functions  $ \rho(r,s) $ and $ \tau(t,s) $ are linear equivalent:

 $$
 C_1 \rho(t,s) \le \tau(t,s) \le C_2 \rho(t,s), \ C_1, C_2 = \const \in (0,\infty). \eqno(3.6)
 $$
Then the entropy condition (3.3) is satisfied. \par
{\bf Proof.} Since the limit process $ \zeta_{\infty}(t)  $ is continuous,  Gaussian  and has at the same covariation function
as $ \eta(t), $ it is stationary in strong sense and also continuous. Our assertion follows from the famous necessary condition
for sample part continuity of Gaussian stationary  process belonging to X.Fernique \cite{Fernique1}. \par
\vspace{3mm}

{\bf Remark 3.1.} At the same result is true even without assumption of stationarity with at the same proof if the function
$  t \to \tau(t,0) $ is monotonic in some neighborhood of zero,  see \cite{Marcus1}. \par

\vspace{4mm}

\section{Some examples.}

\vspace{4mm}

{\bf 1.} Let here $ T = [0, e^{-4}] $ and  let $ \delta = \const \in (0, 1/4).  $ We define the r.p. $ \eta_0(t), \ t \in T $
as follows.

$$
\eta_0(t) := \frac{w(t)}{(2 t)^{1/2} \ ( \log |\log t|)^{ 1/2 + \delta/2} }, \ 0 < t \le e^{-4}, \eqno(4.1)
$$
and $  \eta_0(0): = 0; \ w(t) $  is ordinary Brownian motion. \par
 It follows from the classical Law of Iterated Logarithm (LIL) for Wiener process that
the r.p. $ \eta_0(t) $ is continuous a.e. Since it is "per se" Gaussian, it satisfies the CLT(C(T)). The conditions
(3.1) and (3.2) are also satisfied, but the entropy integral (3.3) is divergent.\par

  Indeed, we have denoting

$$
\tau_0(t,s) = \sqrt{ \Var(\eta_0(t) - \eta_0(s)) }:
$$

$$
\tau_0(t,0) \sim \frac{C}{[\log |\log t|]^{ 1 + \delta  } },
$$
therefore

$$
h_+(\epsilon) := \sup_{t \in T} \mu(B(t,\epsilon)) \ge \mu(B(0,\epsilon)) \ge
C \cdot \exp \left( \exp \left( C_1 \epsilon^{1/(1 + \delta)}  \right) \right), \ 0 < \epsilon < 1/8, \eqno(4.2)
$$
where

$$
B(t,\epsilon) = \{s: s \in T, \ \tau_0(t,s) \le \epsilon \}
$$
and $ \mu  $ is ordinary Lebesgue measure on the real axis.\par
 We can apply the following inequality, see  \cite{Ostrovsky1}, chapter 3,  section 3.2:

$$
\exp H(T,\tau_0, \epsilon) \ge \frac{\mu(T)}{h_+(\epsilon) },
$$
or equally

$$
 H(T,\tau_0, \epsilon) \ge \exp \left( C_2 \epsilon^{1/(1 + \delta) } \right). \eqno(4.3)
$$
 It is easy to verify that for such a entropy the condition (3.3) is not satisfied.\par

\vspace{3mm}

{\bf 2.} The following example appears in the Monte-Carlo method for computation and error estimation of
multiple multivariate parametric integrals, see \cite{Frolov1}, \cite{Grigorjeva1},  \cite{Dudley1}, \cite{Ostrovsky1},
\cite{Ostrovsky106} etc. Namely, let us consider the problem of computation of the following parametric integral

$$
I(t) = \int_D v(t,x) \ \mu(dx),
$$
where $ \mu $ is probability measure:  $ \mu(D) = 1,  $ by means of Monte-Carlo method:

$$
I(t) \approx I_n(t) := n^{-1} \sum_{j=1}^n v(t, \beta_j),
$$
where $  \{\beta_j \}  $ are independent r.v. with distribution $  \mu. $ \par
 Consider for error estimation the following variable:

 $$
 \gamma_n(u) := {\bf P} \left( \sqrt{n} \sup_t | I_n(t) - I(t)  | > u  \right).
 $$
We put here  $ \eta_I(t) = v(t,\beta_1) - I(t).  $  If $ \eta_I(t) $ satisfies the CLT(C(T)), then as $ n \to \infty $

$$
\gamma_n(u) \to  {\bf P} \left( \max_{t \in T} |\zeta_{\infty}(t)| > u \right),
$$
therefore

$$
\gamma_n(u) \approx {\bf P} \left( \max_{t \in T} |\zeta_{\infty}(t)| > u \right) =:\gamma_{\infty}(u).
$$

 The asymptotical behavior or exact estimations as $  u \to \infty $  for the right-hand side of the last equality
are well known, see \cite{Ostrovsky1},chapter 3;  \cite{Piterbarg1}. As a rule as $ u \to \infty $

$$
\gamma_{\infty}(u) \sim K \cdot u^{\kappa-1} \cdot \exp \left( - \frac{u^2}{2 \sigma^2} \right),
$$
where
$$
\sigma^2 = \max_{t \in T} \Var \eta(t) = \max_{t \in T} R(t,t), \ K, \kappa = \const \in (0,\infty).
$$
 Let $ \varepsilon $ be some "small" number, for instance, 0.05 or 0.01 etc.   Denote by $ U(\varepsilon) $ the
maximal root of equation

$$
\gamma_{\infty}(U(\varepsilon)) = \varepsilon,
$$
we deduce that  if $ \eta(t) $ satisfies CLT(C(T)), then with probability $ \approx 1 - \varepsilon $

$$
\sup_{t \in T} |I_n(t) - I(t)| \le \frac{U(\varepsilon)}{\sqrt{n}};
$$
which is a twice asymptotical: as $ n \to \infty $ and as  $ u \to \infty  $ confidence region for $ I(t) $ in the uniform norm.  \par
 The non-asymptotical exponential exact estimation for $ \gamma_{\infty}(u) $ in the modern terms of majorizing measures see in the
article  \cite{Ostrovsky102}.\par
\vspace{3mm}

 Suppose there exist a r.v. $ \theta $ and a non-random  monotonically decreasing sequence $ \delta(n), $ such that

$$
E(n,g) \le \theta \cdot \delta(n), \ \lim_{n \to \infty} \delta(n) = 0. \eqno(4.4)
$$
  We impose the following condition on the r.v. $ \theta: $

$$
\exists m = \const > 1, \ C = \const \in (0,\infty),  \ {\bf P} (\theta > x) \le e^{ -C x^{m}}, \ x \ge 0, \eqno(4.5)
$$
or equally

$$
|\theta|_p \le C_2 \ p^{1/m}, \ p \ge 1,
$$
and denote $ \tilde{m} = \min(m.2), \ m' = \tilde{m}/(\tilde{m}-1).  $  The conditions of theorem 2.2 are satisfied if for example

$$
\sum_{r=1}^{\infty} \frac{\delta(2^r)}{r^{1/\tilde{m}}} < \infty. \eqno(4.6).
$$
 The condition (4.6) is satisfied in turn if for instance for some positive value $ \Delta = \const > 0  $

$$
\delta(n) \le \frac{C_3}{ [\log (n+2)]^{1/m' + \Delta} }. \eqno(4.7)
$$
 We used the following proposition, see \cite{Kozachenko1}, \cite{Ostrovsky1}, chapter 2, section 2.1, page 50-53:
if $  \xi_i, \ i = 1,2, \ldots  $ are centered, i., i.d. r.v.  such that

$$
{\bf P}(|\xi_i| > x) \le \exp \left(  - x^m \right), \ x > 0,
$$
 then

 $$
 \sup_n {\bf P} \left( n^{-1/2} \left| \sum_{i=1}^n \xi_i  \right| > x    \right)   \le \exp \left(  - C(m) \ x^{\tilde{m}} \right), \ x > 0,
 $$
 and the last  inequality is exact.\par

  \vspace{3mm}
  The examples of a (random) functions $ g(\cdot) $ for which

 $$
\delta(n) \asymp \frac{C_4}{ [\log (n+2)]^{1/m' + \Delta} }
$$
 may be constructed by means of lacunary random series, see \cite{Marcus3},  \cite{Ostrovsky301}. \par

 Notice that the conditions (4.6) and (4.7) are so weak so that the entropy  condition \cite{Marcus2}
and majorizing measures conditions \cite{Heinkel1}  are not satisfied. \par

\vspace{3mm}

{\bf Remark 4.1.} The condition (4.4) is satisfied if for example

$$
\omega(g,\epsilon) := \sup_{h: |h| \le \epsilon} || g (\cdot + h) - g(\cdot)  ||   \le C \cdot \theta \cdot
\delta( \Ent[1/\epsilon]),
$$
 Jackson's inequality; see \cite{Achieser1}, \cite{DeVore1},  \cite{Lorentz1}, \cite{Natanson1};
but the inverse inequality is non true, see aforementioned Stechkin's  estimate (\cite{Timan1}, chapter 6, section 6.1). \par

\vspace{3mm}

{\bf 3.} Let now the set $ T = \{1,2,\ldots; \infty \}   $ equipped with the distance

$$
d(i,j) = \left| \frac{1}{j}  -  \frac{1}{i}   \right|, \ d(i, \infty) = \frac{1}{i}. \ d(\infty,\infty) = 0.\eqno(4.8)
$$
 The set $ T $ is compact set relatively  the distance $ d. $ \par
  We choose in the capacity of probability space the ordinary interval $ (0,1) $ with Lebesgue measure. \par
  Let $ p_0 = \const \in (1,2), $

$$
 \alpha = \const \in (0, \min(1, p_0/(2-p_0)), \  a(n) = 1 - 0.5 n^{-\alpha}, \ n = 1,2,\ldots;  \
$$

$$
\Delta(n) = a(n+1) - a(n), \ c(n) = n^{\alpha/p_0}.
$$

Let also

$$
f_{1/2}(x) = \sqrt{|\log x|}, \ x \in (0,1); \ f_{1/2}(x) = 0, \ x \notin (0,1). \eqno(4.9)
$$
 We define for the values $  n \in T $ the following r.p.

$$
\eta_n(x) = c(n) \ \epsilon_n \ f_{1/2} \left(\frac{x-a(n)}{\Delta(n)} \right),  \ \eta_{\infty}(x) = 0, \eqno(4.10)
$$
where $ \{ \epsilon_n \} $ is a Rademacher's sequence defined on some other probabilistic space  and independent
on $ f_{1/2} (\cdot), $ so that $ {\bf E} \eta_n = 0.  $\par
 Recall that

 $$
 ||\eta(\cdot)|| = \sup_n |\eta_n|.
 $$

 This example was offered by the authors in \cite{Ostrovsky104} for another purpose. It was proved in particular
in \cite{Ostrovsky104}  that the r.p. $ \eta_n $ is continuous with probability one, this imply:

$$
{\bf P} \left( \lim_{n \to \infty} \eta_n \to 0 \right) = 1,
$$
and
$$
{\bf E} ||\eta||^{p_0} < \infty, \ \forall p > p_0 \ \Rightarrow {\bf E} ||\eta||^{p} = \infty.  \eqno(4.11)
$$

 Let us prove in addition that the sequence of a r.v. $ \{ \eta_n \}  $ is pre-gaussian in the considered space $ C(T,d). $
This imply by definition that the mean zero Gaussian distributed sequence  $ \{ \nu_n \}  $
with at the same covariation function $ Q(n_1, n_2)  $ as $ \{ \eta_n \} $ belong to the set $ C(T,d) $ with probability one.\par
 Namely, since the Rademacher's sequence contains from independent r.v. $ \epsilon_n, $

$$
Q(n_1,n_2) = 0, \ n_1 \ne n_2.
$$
 It suffices  to prove $ {\bf P} (\nu_n \to 0) = 1. $ We have:

 $$
 \Var{\nu_n} = Q(n,n) = \Var{\eta_n} = {\bf E} |\eta_n|^2.
 $$
It is proved in \cite{Ostrovsky104} that

$$
{\bf E} |\eta_n|^2 \le C_1 \ n^{-C_2}, \ C_2 > 0,
$$
therefore $  \Var{\nu_n} \le C_1 \ n^{-C_2} $ and by virtue of Gassiness of the sequence $ \{ \nu_n  \}  $

$$
\forall \epsilon > 0 \ \Rightarrow \  \sum_{n=1}^{\infty} {\bf P} ( | \nu_n| > \epsilon  )< \infty.
$$
Thus,

$$
{\bf P} \left( \lim_{n \to \infty} \nu_n \to 0 \right) = 1.
$$

 The r.p. $ \eta_n $ does not satisfy the majorizing measure condition, \cite{Ostrovsky104}. It follows from
(4.11) that  the condition (3.4) is also not satisfied. Therefore, $ \eta_n $ does not satisfy CLT in our space
$ C(T.d). $\par

\vspace{4mm}

\section{Concluding remarks.}

\vspace{4mm}

{\bf 1.} Note that the r.v. $ M $ and the distance $ \rho = \rho(t,s) $ in (3.1) may be introduced constructively in two stages
as follows. First step: define the r.v. $ L $

$$
L := \sup_{t,s \in T} |\eta(t) - \eta(s)|.
$$
 Second step:

 $$
 q(t,s) := \vraisup_{\Omega} \left[ \frac{|\eta(t) - \eta(s)|}{L} \right];
 $$
(the case $ L=0 $ is trivial). Then the function $ (t,s) \to q(t,s)  $  is continuous bounded distance
(more exactly, semi-distance) on the set $  T $  and

$$
|\eta(t) - \eta(s)| \le L \cdot q(t,s).
$$
 Obviously, this choice of the variables  $ L, \ q  $ is optimal.\par

\vspace{4mm}

{\bf 2.}  It is clear that we can use instead trigonometrical approximation  the approximation  by means of
algebraic polynomials.\par

\vspace{4mm}

{\bf 3.}  Suppose that  instead the condition (3.6) is true the following inequality:

$$
\rho(t,s) \le C \cdot \tau^{1/\beta} (t,s),      \ \beta = \const \in (0,1).
$$
 Then the {\it necessary}  condition of a view (3.3) must be transformed as follows:

$$
\int_0^1 H^{1/2}( T,\rho, z^{\beta}) \ dz < \infty.
$$

\vspace{4mm}

{\bf 4.} It is no hard to generalize our results into the multidimensional case $ T = [0, 2 \pi]^d  $
and into the non-periodical continuous function space.\\


\vspace{4mm}

\begin{thebibliography}{99}

\vspace{4mm}

\bibitem{Achieser1}
{\sc Achieser N.I.} {\it Theory of Approximation.}
 Frederic Ungar Publishing Co., Second Printing, New York, (1956).

\bibitem{Bednorz1}
 {\sc Bednorz W.} (2006). {\it A theorem on Majorizing Measures.}
 Ann. Probab., {\bf 34}, 1771-1781. MR1825156.

\bibitem{Bednorz2}
{\sc Bednorz W.} {\it The majorizing measure approach to the sample boundedness.}
arXiv:1211.3898v1 [math.PR] 16 Nov 2012

\bibitem{Bednorz3}
{\sc Bednorz W.}  (2010), {\it Majorizing measures on metric spaces. } C.R. math. Acad. Sci.
Paris, (2010), 348, no. 1-2, 75-78, MR2586748

\bibitem{Bednorz4}
{\sc Bednorz W.}  {\it H\"older continuity of random processes.}
arXiv:math/0703545v1 [math.PR] 19 Mar 2007.

\bibitem{Belyaev1}
{\sc Belyaev Yu.K.} {\it  Continuity and H\"older condition  for sample functions of stationary Gaussian
processes.  } 4th Berkeley Symp. Stat. Probab.  Vol. 2 (1961) p. 23-33.

\bibitem{Buldygin1}
{\sc Buldygin V.V.} {\it  Supports of probability measures in separable Banach Spaces. } Probab. Theory
Appl., (1984), V. 29 Issue 3, p. 528-532 (in Russian).

\bibitem{Buldygin2}
{\sc Buldygin V.V.} {\it Convergence of random elements in topological spaces.}  Naukova Dumka, Kiev, (1980),
(in Russian).

\bibitem{Poussin1}
 {\sc Ch.J. de la Vallee-Poussin}. {\it Sur la meilleure approximation des fonctions d'une variable reelle par
 des expressions d'ordre donne.} C.R. Acad. Sci. Paris Sér. I. Math. , 166 (1918) pp. 799–802

\bibitem{DeVore1}
{\sc De Vore A., Lorentz G.G.} {\it Constructive approximation.}  Springer, (1993), Berlin, Heidelberg, New York.

\bibitem{Dudley1}
{\sc Dudley R.M.} {\it Uniform Central Limit Theorem.}  Cambridge University Press, 1999.

\bibitem{Fernique1}
 {\sc Fernique X.} (1975). {\it Regularite des trajectoires des
    function aleatiores gaussiennes.}  Ecole de Probablite de
    Saint-Flour, IV - 1974, Lecture Notes in Mathematic. {\bf 480}, 1 - 96,
    Springer Verlag, Berlin.

\bibitem{Fernique2}
{\sc Fernique X,} {\it Caracterisation de processus de trajectoires majores ou
continues.} Seminaire de Probabilit´s XII. Lecture Notes in Math. 649,
(1978), 691–706, Springer, Berlin.

\bibitem{Fernique3}
{\sc Fernique X.} {\it Regularite de fonctions aleatoires non gaussiennes.} Ecolee
de Ete de Probabilit´s de Saint-Flour XI-1981. Lecture Notes in Math.
976, (1983), 1–74, Springer, Berlin.

\bibitem{Frolov1}
{\sc Frolov A.S., Tchentzov N.N.} {\it On the calculation by the Monte-Carlo
method definite integrals depending on the parameters.} Journal of Computetional
Mathematics and Mathematical Physics, (1962), V. 2, Issue 4, p. 714-718
(in Russian).

\bibitem{Garsia1}
{\sc Garsia, A. M.; Rodemich, E.; and Rumsey, H., Jr.} {\it A real variable lemma and
the continuity of paths of some Gaussian processes.} Indiana Univ. Math. J. 20
(1970/1971), 565-578.

\bibitem{Gikhman1}
{\sc Gikhman I.I., Skorohod A.V. } {\it  Introduction to the theory of random processes. }
Nauka, GIFML, Moscow (1965).

\bibitem{Gine1}
{\sc Gin\'e E.} {\it  On the Central Limit theorem for sample continuous processes. }
Ann. Probab. (1974), {\bf 2,} 629-641.

\bibitem{Grigorjeva1}
{\sc Grigorjeva M.L., Ostrovsky E.I.}  {\it Calculation of Integrals on discontinuous
Functions by means of depending trials method.} Journal of Computational
Mathematics and Mathematical Physics, (1996), V. 36, Issue 12, p. 28-39 (in
Russian).

\bibitem{Heinkel1}
{\sc Heinkel B.} {\it  Measures majorantes et le theoreme de la limite centrale dans $ C(S). $ }
 Z. Wahrscheinlichkeitstheory. verw. Geb., (1977). {\bf 38}, 339-351.

\bibitem{Ibragimov1}
{\sc Ibragimov I.A.}  {\it Properties of sample functions of stochastic processes and embedding
theorems.}  Teory Probab. Appl., 18:3 (1973), 468–480.

\bibitem{Marcus2}
{\sc Jain N.C. and Marcus M.B.} {\it   Central limit theorem for $ C(S) $ valued random variables. }
J. of Funct. Anal., (1975), {\bf 19}, 216-231.

\bibitem{Kozachenko1}
 {\sc Kozachenko Yu. V., Ostrovsky E.I.} (1985). {\it The Banach Spaces of
      random Variables of subgaussian type.} Theory of Probab. and Math.
      Stat. (in Russian). Kiev , KSU, {\bf 32}, 43 - 57.

\bibitem{Kwapien1}
{\sc Kwapien S. and Rosinsky J.} {\it Sample H\"older continuity of stochastic processes and majorizing measures.}
 (2004).  Seminar on Stochastic Analysis, Random Fields
and Applications IV, Progr. in Probab. 58, 155–163. Birkh\"ouser, Basel.

\bibitem{Ledoux1}
 {\sc Ledoux M., Talagrand M.} (1991) {\it Probability in Banach Spaces.}
      Springer, Berlin, MR 1102015.

\bibitem{Lorentz1}
{\sc Lorentz G.G.} {\it Approximation of a Functions.} Chelsea,  (1986), London, Berlin, New York, Hong Kong.

\bibitem{Marcus1}
{\sc Marcus M.B. and Shepp L.A.} {\it  Continuity of Gaussian processes.} Trans. Amer. Math. Soc.,
{\bf 151}, (1970), 377-391.

\bibitem{Marcus3}
{\sc Marcus M.B.} {\it  Gaussian lacunary series and the modulus of continuity for Gaussian processes.}
 Zeitschrift Wahr. Theory und verw. Gebiete, {\bf 22}, 3, (1972), 301-322.

\bibitem{Natanson1}
{\sc Natanson I.P.} {\it Constructive Function Theory.}
 Frederic Ungar Publishing Co.,  New York, (1964).

\bibitem{Ostrovsky1}
{\sc  Ostrovsky E.I.} (1999). {\it Exponential estimations for random Fields and its
applications (in Russian).}  Moscow - Obninsk, OINPE.

\bibitem{Ostrovsky100}
{\sc  Ostrovsky E. and Sirota L.} {\it Module of continuity for the functions
 belonging to the Sobolev-Grand Lebesgue Spaces.}
arXiv:1006.4177v1 [math.FA] 21 Jun 2010

\bibitem{Ostrovsky101}
 {\sc Ostrovsky E., Sirota L.}
{\it Continuity of Functions belonging to the fractional Order Sobolev's-Grand Lebesgue Spaces.}
arXiv:1301.0132v1 [math.FA] 1 Jan 2013

\bibitem{Ostrovsky102}
{\sc Ostrovsky E., Rogover E.} {\it Exact exponential Bounds for the random Field Maximum Distribution
via the Majorizing Measures (Generic Chaining.)}
arXiv:0802.0349v1 [math.PR] 4 Feb 2008

\bibitem{Ostrovsky103}
{\sc Ostrovsky E.I.} (2002). {\it Exact exponential Estimations for Random Field
Maximum Distribution.} Theory Probab. Appl. 45 v.3, 281 - 286.

\bibitem{Ostrovsky104}
 {\sc Ostrovsky E., Sirota L.} {\it A counterexample to a hypothesis of light tail of maximum distribution
 for continuous random processes with light finite-dimensional tails.  }
arXiv:1208.6281v1 [math.PR] 30 Aug 2012

\bibitem{Ostrovsky105}
{\sc Ostrovsky E., Rogover E.} {\it Non - asymptotic exponential bounds for
MLE deviation under minimal conditions via classical and generic chaining methods.}
arXiv:0903.4062v1 [math.PR] 24 Mar 2009

\bibitem{Ostrovsky106}
 {\sc Ostrovsky E., Sirota L.}
{\it Monte-Carlo method for multiple parametric integrals
calculation and solving of linear integral Fredholm equations
of a second kind, with confidence regions in uniform norm.}
arXiv:1101.5381v1 [math.FA] 27 Jan 2011

\bibitem{Ostrovsky207}
 {\sc Ostrovsky E., Sirota L.}
{\it Simplification of the majorizing method,  with development.}
arXiv:1302.3202v1 [math.PR] 13 Feb 2013

\bibitem{Ostrovsky208}
{\sc Ostrovsky E.} {\it  On the local structure of normal fields. }
Soviet Math. Doklady, (1970), v. 105 No 1 p. 1425-1428.

\bibitem{Ostrovsky209}
{\sc Ostrovsky E.} {\it  Convergence of canonical  expression for normal fields. }
 Math Notes, (1973),  V. 14 Issue 4 p. 565-572.

\bibitem{Ostrovsky210}
{\sc Ostrovsky E.} {\it Local properties of Gaussian processes and fields.}
Ph. D. Dissertation, MSU, (1972).

\bibitem{Ostrovsky211}
{\sc Ostrovsky E.} {\it  Covariation operators and some estimated of Gaussian vectors in
Banach space.}
Soviet Math. Doklady, (1977), v. 236 No 3 p. 541-544.

\bibitem{Ostrovsky301}
 {\sc Ostrovsky E., Sirota L.} {\it Theory of approximation and continuity of random processes.    }
arXiv:1303.3029v1 [math.PR] 12 Mar 2013.

\bibitem{Ostrovsky302}
{\sc Ostrovsky E.} {\it  About the support of probabilistic measures
in separable Banach spaces. } Soviet Math. Doklady, (1980), v. 255 No 6 p. 836-838.

\bibitem{Piterbarg1}
{\sc Piterbarg V.I.} {\it Asymptitical methods in the theory of Gaussian processes and fields.  }
Moscow, MSU, (1988), (in Russian).

\bibitem{Prokhorov1}
{\sc Prokhorov Yu.V.}   {\it Convergense of Random Processes and Limit Theorems
of Probability Theory.} Probab. Theory Appl., (1956), V. 1, 177-238.

\bibitem{Ral'chenko1}
{\sc Ral'chenko, K. V. } {\it The two-parameter Garsia-Rodemich-Rumsey inequality and its
application to fractional Brownian fields.} Theory Probab. Math. Statist. No. 75
(2007), 167-178.

\bibitem{Rivlin1}
{\sc Rivlin T.J.} {\it Introduction to the Approxination of Functions. } Dover, (1969), New York.

\bibitem{Talagrand1}
 {\sc Talagrand M.} (1996). {\it Majorizing measure: The generic chaining.}
 Ann. Probab., {\bf 24} 1049 - 1103. MR1825156

\bibitem{Talagrand2}
 {\sc Talagrand M.} (2005). {\it The Generic Chaining. Upper and
     Lower Bounds of Stochastic Processes.} Springer, Berlin. MR2133757.

\bibitem{Talagrand3}
{\sc Talagrand M.} (1987).{\it Regularity of Gaussian processes.} Acta Math. 159 no. 1-2, 99–
149, MR 0906527.

\bibitem{Talagrand4}
{\sc Talagrand M.} (1990),  {\it Sample boundedness of stochastic processes under increment
conditions. } Annals of Probability 18, N. 1, 1-49, MR1043935.

\bibitem{Talagrand5}
{\sc Talagrand M.} (1992). {\it A simple proof of the majorizing measure theorem. } Geom.
Funct. Anal. 2, no. 1, 118 \ - \ 125. MR 1143666

\bibitem{Timan1}
{\sc Timan A.F.} {\it Theory of Approximation of a real variables.}  Pergamon Press, Oxford, London, New York, Paris, (1963)

\bibitem{Watanabe1}
{\sc Watanabe H.} {\it On the continuity property of Gaussian random fields. } Studia Mathematica, V. XL1X,
(1973), 81-90.

\end{thebibliography}
\end{document}